\documentclass[12pt,twoside]{article}
\usepackage[english]{babel}
\usepackage{amsmath}
\usepackage{times,amssymb,amscd}
\usepackage[utf8]{inputenc}
\usepackage[pass,legalpaper]{geometry}
\usepackage{fancyhdr}
\usepackage{multicol}
\usepackage{enumitem}
\usepackage{graphicx}
\usepackage{pgf}
\usepackage{textcomp}
\newtheorem{thm}{Theorem}[section]

\newtheorem{rem}[thm]{Remark}
\numberwithin{equation}{section}

\newcommand{\bH}{\mathbb{H}}

\newcommand{\bs}{\mathbf{s}}
\newcommand{\ba}{\mathbf{a}}

\newcommand{\bx}{\mathbf{x}}

\newcommand{\bu}{\mathbf{u}}

\newcommand{\Bu}{\boldsymbol{u}}

\newcommand{\Bs}{\boldsymbol{s}}

\newcommand{\HYP}{\bH^3}

\begin{document}

	\title{Visualization of sphere and horosphere packings related to Coxeter tillings generated by simply truncated orthoschemes with parallel faces}
	\author{Arnasli Yahya, Jen\H{o} Szirmai \\ 
\normalsize Budapest University of Technology and Economics, Institute of Mathematics, \\
\normalsize Department of Geometry \\
\normalsize Budapest, P. O. Box: 91, H-1521 }
	\maketitle

\begin{abstract}

	In this paper we describe and visualize the densest ball and horoball packing configurations belonging to 
	the simply truncated $3$-dimensional hyperbolic Coxeter orthoschemes with parallel faces using the results of \cite{YSz}. 
	These beautiful packing arrangements describe and show the very interesting structure of the mentioned orthoschemes and the corresponding Coxeter groups.
	We use for the visualization the Beltrami-Cayley-Klein ball model of $3$-dimensional hyperbolic space $\HYP$ and the pictures were made by the Python software.
	
\end{abstract}
\section{Introduction}
Visualization of mathematical problems is not only a representation of specific objects or an approach in the teaching process, 
but also plays an important role in understanding the problem and developing solution steps. 
It can be shown the deeper contexts of the problem and the possibilities for moving forward. 

In hyperbolic spaces $\mathbb{H}^n$ for $2 \leq n \leq 9$, the known densest ball and horoball configurations are derived by Coxeter simplex tilings, 
therefore it is interesting to determine their optimal horoball packings related to Coxeter tilings. In the former papers, we investigated that Coxeter simplex tilings whose generating simplices do not have parallel faces. In \cite{YSz} we extended our study to the $3$ dimensional Coxeter tilings generated by simple frustum orthoschemes with parallel faces. But first, we recall some results related to the above topic.

In the case of periodic ball or horoball packings, this local density defined above can be extended to the entire hyperbolic space and 
is related to the simplicial density function (defined below) that we generalized in \cite{SzJ2} and \cite{SzJ3}. In this paper, we shall use such definition of packing density.  

A Coxeter simplex is a top dimensional simplex in $\overline{\mathbb{H}}^n$ with dihedral angles either integral submultiples of $\pi$ or zero. 
The group generated by reflections on the sides of a Coxeter simplex is a Coxeter simplex reflection group. 
Such reflections generate a discrete group of isometries of $\mathbb{H}^n$ with the Coxeter simplex as the fundamental domain; 
hence the groups give regular tessellations if the fundamental simplex is characteristic. We note here that the Coxeter groups 
are finite for $\mathbb{S}^n$, and infinite for $\mathbb{E}^n$ or $\overline{\mathbb{H}}^n$. 

There are non-compact Coxeter simplices in $\overline{\mathbb{H}}^n$ with ideal vertices in $\partial \mathbb{H}^n$, 
however only for dimensions $2 \leq n \leq 9$; furthermore, only a finite number exists in dimensions $n \geq 3$. 
Johnson {\it et al.} \cite{JKRT} found the volumes of all mentioned Coxeter simplices in hyperbolic $n$-space, also see Kellerhals \cite{KH1}. 
Such simplices are the most elementary building blocks of hyperbolic manifolds, the volume of which is an important topological invariant. 

In $n$-dimensional space of constant curvature $(n \geq 2)$, define the simplicial density function $d_n(r)$ to be the density of $n+1$ mutually tangent balls of radius $r$ in the simplex spanned by their centers. L.~Fejes T\'oth and H.~S.~M.~Coxeter conjectured that the packing density of balls of radius $r$ in $X$ cannot exceed $d_n(r)$. Rogers \cite{Ro64} proved this conjecture in Euclidean space $\mathbb{E}^n$. The $2$-dimensional spherical case was settled by L.~Fejes T\'oth \cite{FTL}, and B\"or\"oczky \cite{BK3}, who extend the analogous statement to  $n$-dimensional spaces of constant curvature.

The simplicial packing density upper bound $d_3(\infty) = (1+\frac{1}{2^2}-\frac{1}{4^2}-\frac{1}{5^2}+\frac{1}{7^2}+\frac{1}{8^2}--++\dots)^{-1} = 0.85327\dots$ cannot be achieved by packing regular balls, instead it is realized by horoball packings of $\overline{\mathbb{H}}^3$, the regular ideal simplex tiles $\overline{\mathbb{H}}^3$. More precisely, the centers of horoballs in  $\partial\overline{\mathbb{H}}^3$ lie at the vertices of the ideal regular Coxeter simplex tiling with Schl\"afli symbol $\{3,3,6\}$. 

In \cite{KSz} we proved that this optimal horoball packing configuration in $\mathbb{H}^3$  is not unique. We gave several more examples of regular horoball packing arrangements based on asymptotic Coxeter tilings using horoballs of different types, that is horoballs that have different relative densities with respect to the fundamental domain, that yield the B\"or\"oczky--Florian-type simplicial upper bound \cite{BK4}.
 
Furthermore, in \cite{SzJ2,SzJ3} we found that 
by allowing horoballs of different types at each vertex of a totally asymptotic simplex and generalizing 
the simplicial density function to $\overline{\mathbb{H}}^n$ for $(n \ge 2)$, the B\"or\"oczky-type density upper bound is not valid for the fully asymptotic simplices for $n \geq 4$. 
For example, in $\overline{\mathbb{H}}^4$ the locally optimal simplicial packing density is $0.77038\dots$, higher than the B\"or\"oczky-type density upper bound of $d_4(\infty) = 0.73046\dots$ using horoballs of a single type. 
However these ball packing configurations are only locally optimal and cannot be extended to the entirety of the
ambient space $\overline{\mathbb{H}}^n$. 
In \cite{KSz2} we found seven horoball packings of Coxeter simplex tilings in $\overline{\mathbb{H}}^4$ that 
yield densities of $5\sqrt{2}/\pi^2 \approx 0.71645$, counterexamples to L. Fejes T\'oth's conjecture of $\frac{5-\sqrt{5}}{4}$ stated in his foundational 
book {\it Regular Figures} \cite[p. 323]{FTL}.

In \cite{KSz3}, \cite{KSz4} we extend our study of horoball packings to $\overline{\mathbb{H}}^n$ ($5\le n \le 9$) using our methods that were successfully applied in lower dimensions. In the previously mentioned papers we studied the ball and horoball packing 
related to the Coxeter simplex tilings where the vertices of the simplices 
are proper points of the hyperbolic space 
$\mathbb{H}^n$ or they are ideal i.e. lying on 
the sphere $\partial \mathbb{H}^n$.

In \cite{YSz}, we considered the Coxeter tilings in $3$-dimensional hyperbolic space $\mathbb{H}^3$ where the generating orthoscheme is a simple truncated Coxeter orthoscheme with parallel faces i.e. their dihedral angle is zero. Here we studied the Coxeter tilings which are given with Schl\"afli symbol $\{\infty,q,r,\infty\}$ (see Fig.~1. second graph). We determined their optimal ball and horoball packings, proved that the densest packing arrangement of the considered tilings is realized at the tilings $\{\infty,3,6,\infty\}$, and $\{\infty;6;3;\infty \}$ by horoballs with density $\approx 0.8413392$. 
\begin{rem}
The first cases in Fig.~1. is investigated in papers \cite{SzJ9,SzJ10,SzJ11} where from the truncated orthoschemes can be derived prism like tilings generated hyperball or hyp-hor packings
(or coverings) in $\mathbb{H}^3$. 
\end{rem}
In this paper we describe and visualize the densest ball and horoball packing configurations belonging to 
the simply truncated $3$-dimensional hyperbolic Coxeter orthoschemes with parallel faces using the results of \cite{YSz}. 
These beautiful packing arrangements describe and show the very interesting structure of the mentioned orthoschemes and the corresponding Coxeter groups.
We use for the visualization the Beltrami-Cayley-Klein ball model of $3$-dimensional hyperbolic space $\HYP$ and the pictures were made by the Ipy Volume, a Python library for visualizing 3D-object.
\section{Basic Notions}
For the computations and visualization, we use the projective model of the hyperbolic space $\HYP$. The model is defined in the $\mathbb{E}^{1,n}$ Lorentz space with signature $(1,n)$, i.e. consider $\mathbf{V}^{n+1}$ real vector space equipped with the bilinear form: 
\begin{equation*}
\langle ~ \mathbf{x},~\mathbf{y} \rangle = -x^0y^0+x^1y^1+ \dots + x^n y^n.
\end{equation*}
In the vector space, consider the following equivalence relation:
\begin{equation*}
\mathbf{x}(x^0,...,x^n)\thicksim \mathbf{y}(y^0,...,y^n)\Leftrightarrow \exists~ c\in \mathbb{R}\backslash\{0\}: \mathbf{y}=c\cdot\mathbf{x}.
\end{equation*}
The factorization with $\thicksim$ induces the $\mathcal{P}^n(\mathbf{V}^{n+1},\mbox{\boldmath$V$}\!_{n+1})$ $n$-dimensional real projective space. In this space to interpret the points of $\mathbb{H}^n$ hyperbolic space, consider the following quadratic form:
\begin{equation*}
Q=\{[\mathbf{x}]\in\mathcal{P}^n | \langle ~ \mathbf{x},~\mathbf{x} \rangle =0 \}=:\partial \mathbb{H}^n.
\end{equation*}
The inner points relative to the cone-component determined by $Q$ are the points of $\mathbb{H}^n$ (for them $\langle ~ \mathbf{x},~\mathbf{x} \rangle <0$), 
the point of $Q=\partial \mathbb{H}^n$ are called the points at infinity, and the points lying outside relative to $Q$ are outer points of $\mathbb{H}^n$ 
(for them $\langle ~ \mathbf{x},~\mathbf{x} \rangle >0$). We can also define a linear polarity between the points and hyperplanes of the space: 
the polar hyperplane of a point $[\mathbf{x}]\in\mathcal{P}^n$ is 
$Pol(\mathbf{x}):=\{[\mathbf{y}]\in\mathcal{P}^n | \langle ~ \mathbf{x},~\mathbf{y} \rangle =0 \}$, and hence $\mathbf{x}\in\mathbf{V}^{n+1}$ is 
incident with ${\boldsymbol{a}}\in\mbox{\boldmath$V$}\!_{n+1}$ iff $ \mathbf{x}~{\boldsymbol{a}} =0$. In this projective model, 
we can define a metric structure related to the above bilinear form, where for the distance of two proper points:
\begin{equation}
\cosh\left(\frac{d(\mathbf{x},\mathbf{y})}{k}\right)=\frac{-\langle ~ \mathbf{x},~\mathbf{y} \rangle}{\sqrt{\langle ~ \mathbf{x},~\mathbf{x} \rangle \langle ~ \mathbf{y},~\mathbf{y} \rangle}},
~ (\text{at present}~ k=1).
\end{equation}
This corresponds to the distance formula in the well-known Beltrami-Cayley-Klein model.
\section{The structures of truncated orthoschemes}
Our aim to visualize the simply asymptotic orthoschemes that contain parallel hyperplane faces in 3-dimensional hyperbolic space.
This orthoschemes are represented by their Coxeter graphs (see Fig.\ref{Coxeter graph}),
\begin{figure}[h!]
\begin{centering}
    \includegraphics[scale=0.4]{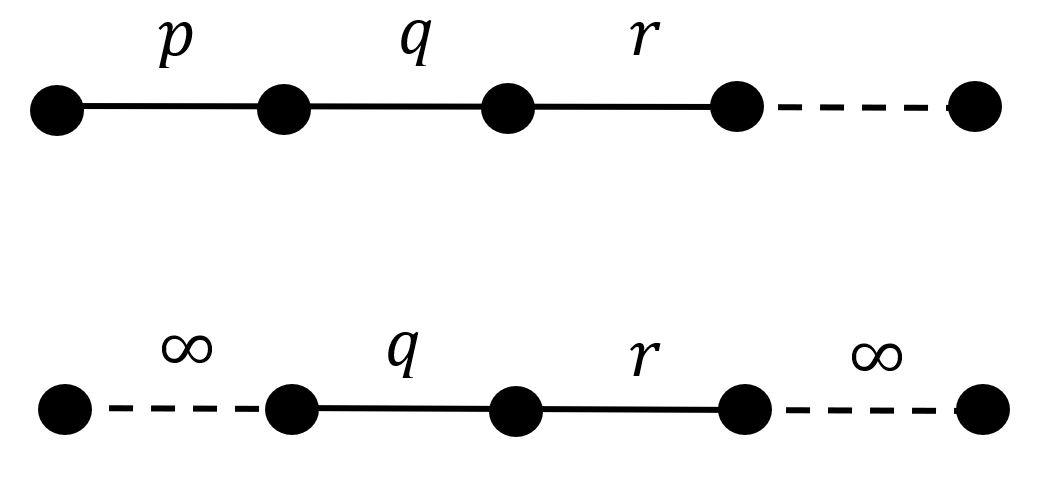}
    \caption{The Coxeter graph}
    \label{Coxeter graphs}
\end{centering}
\end{figure}
where the parameters $p, q, r$ are satisfy the inequalities $\frac{1}{p}+\frac{1}{q}<\frac{1}{2}$ and $\frac{1}{q}+\frac{1}{r} \geq \frac{1}{2}$. 

We will study the second type of the orthoschemes that has corresponding Coxeter-Schl\"{a}fli Matrix as follows
\begin{equation} \label{CoxeterMatrix}
    C=\begin{bmatrix}
    1 & -\cos{(\frac{\pi}{p})}&0&0&0\\
    -\cos{(\frac{\pi}{p})}&1&-\cos{(\frac{\pi}{q})}&0&0\\
    0&-\cos{(\frac{\pi}{q})}&1&-\cos{(\frac{\pi}{r})}&0\\
    0&0&-\cos{(\frac{\pi}{r})}&1&c_4\\
    0&0&0&c_4&1
    \end{bmatrix},
\end{equation}
where the constant $c_4$ can be uniquely determined in the arrangement of Napier cycles \cite{IH2}, i.e
\begin{frame}{}
\begin{equation*}
    c_4 = -\sqrt{\frac{1+\cos^2(\frac{\pi}{p})\cos^2(\frac{\pi}{r})-\cos^2(\frac{\pi}{p})-\cos^2(\frac{\pi}{q})-\cos^2(\frac{\pi}{r})}{1-\cos^2(\frac{\pi}{p})-\cos^2(\frac{\pi}{q})}}.
\end{equation*}
\end{frame}
In our case, there are two parallel faces (hyperplanes) that meet in an ideal point. That means the dihedral angle between these two hyperplanes is equal to $0$.
Therefore, we assume that these two hyperplanes are $h_0$ and $h_1$. Thus, their dihedral angle is $w_0=\frac{\pi}{p} \rightarrow 0$, if $p$ tends to $\infty$, the Coxeter-Schl\"afli matrix (\ref{CoxeterMatrix}) would change to the following form
\begin{equation}\label{CexeterInfinity}
    C=\begin{bmatrix} 
    1 & -1&0&0&0\\
    -1&1&-\cos{(\frac{\pi}{q})}&0&0\\
    0&-\cos{(\frac{\pi}{q})}&1&-\cos{(\frac{\pi}{r})}&0\\
    0&0&-\cos{(\frac{\pi}{r})}&1&-1\\
    0&0&0&-1&1
    \end{bmatrix}
\end{equation}
We notice that $C$ is a singular matrix. The crucial issue is that the complete orthoscheme is uniquely, up to isometry, determined by the 
$4 \times 4$-non singular principal submatrix of $C$, we write this principal submatrix as $C(p,q,r)$,
\begin{equation}\label{CoxeterSubmatrix}
    C(p,q,r)=\begin{bmatrix}
1 & -\cos{(\frac{\pi}{p})}&0&0 \\ -\cos{(\frac{\pi}{p})}&1&-\cos{(\frac{\pi}{q})} & 0\\ 0&-\cos{(\frac{\pi}{q})}&1&-\cos{(\frac{\pi}{r})} \\ 0&0&-\cos{(\frac{\pi}{r})}&1
\end{bmatrix}.
\end{equation}
By the Proposition 1.6 in \cite{IH2}, one could have a set (unique, up to isometries) of 3+1 unit positive vectors 
$\Bu_0, \Bu_1, \Bu_2, \Bu_3$ admitting $C(p,q,r)$ as their Gram matrix and thus generating a Napier cycle in $\mathbb{R}^{1,3}$.
Their pols $\bu_0, \bu_1, \bu_2, \bu_3$ are the ``unit normals" of the hyperplanes $h_0, h_1, h_2, h_3$ respectively,
\begin{equation*}
    h_i=\{\bx\in \mathbb{R}^{1,3} |\langle \bx, \bu_i \rangle =0\},
\end{equation*}
for $i=0,1,2,3$.
The matrix $A$, inverse of matrix $C(p,q,r)$, would uniquely determine the coordinates of the orthoscheme vertices relative to chosen basis vectors $\bu_0, \bu_1,\bu_2, \bu_3$.
The location of the vertices is completely described by the column of matrix $EA$, where $E=\left[\bu_0,\bu_1,\bu_2,\bu_3\right]$, therefore by this observation, 
we choose the initial vectors $\bu_0=\begin{bmatrix}
    \sinh{t} \\ 0 \\ \cosh{t} \\ 0
    \end{bmatrix}$ and $\bu_3=\begin{bmatrix}
    0 \\ 0 \\ 0 \\ 1
    \end{bmatrix}$
where the parameter $t \in \mathbb{R}$. We determine the vectors $\bu_1$, $\bu_2$, and $\bu_4$ based on these initial vectors and by the benchmark Coxeter-Schl\"afli matrix $C$.

\begin{figure}[h!]
\begin{center}
  \begin{minipage}[b]{0.4\textwidth}
  \centering
   \includegraphics[scale=0.7]{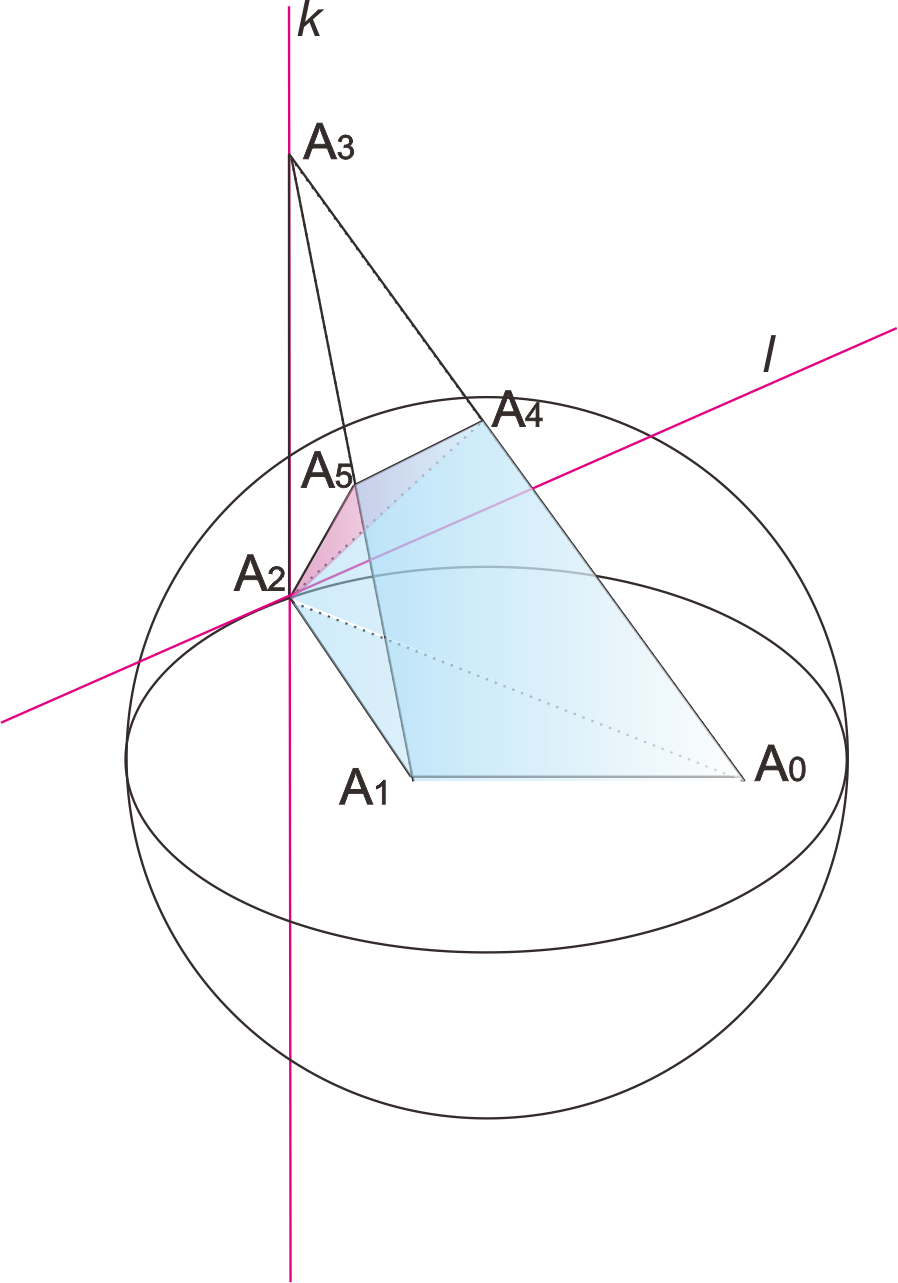}
    \caption{The structure of truncated orthoscheme with the two intersection pairs of its parallel faces.}
    \label{Structure_Orthoscheme}
    \end{minipage}~~~~
    \begin{minipage}[b]{0.4\textwidth}
  \centering
   \includegraphics[scale=0.7]{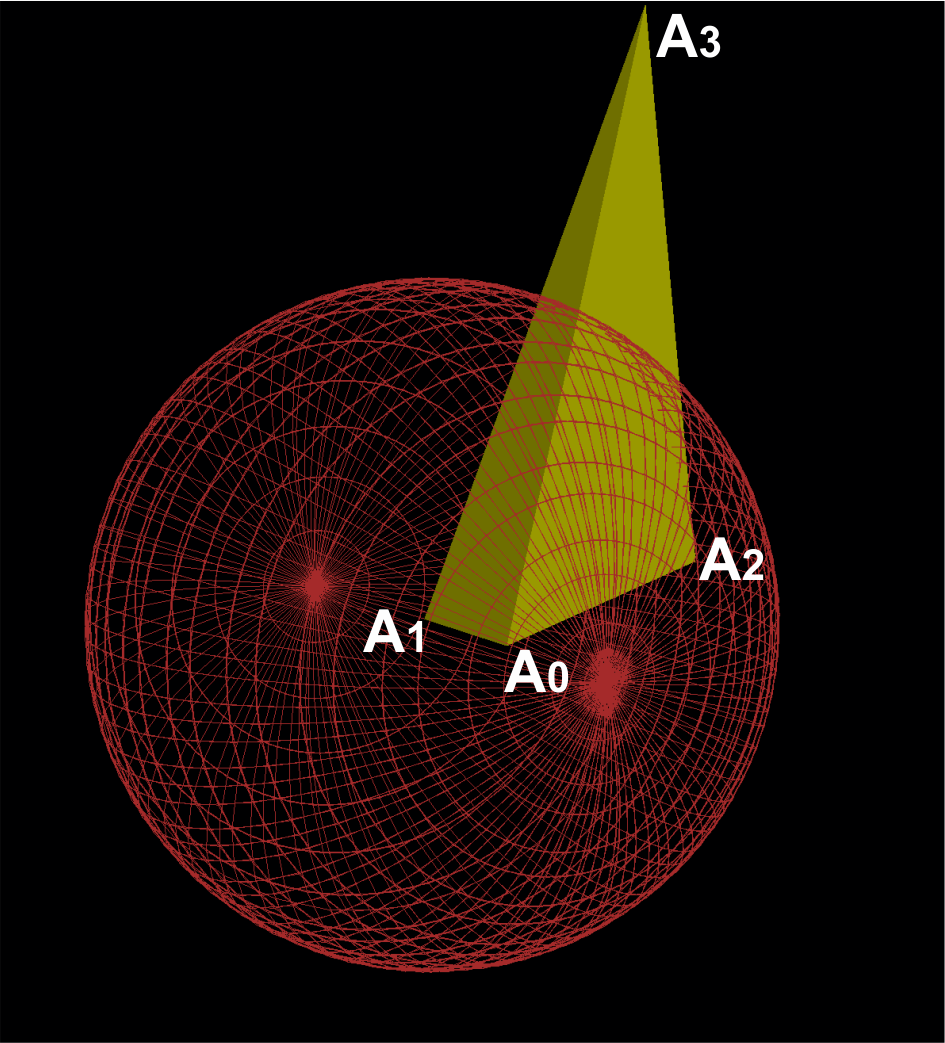}
    \caption{Orthoscheme with a ultra ideal vertex $A_3$, and the couple of parallel faces ($A_0 A_2 A_3$ and $A_2 A_1 A_3$)}
    \label{Orthoscheme}
    \end{minipage}
    \end{center}
\end{figure}

\begin{figure}
    \begin{center}
        \begin{minipage}[b]{0.4\textwidth}
            \centering
            \includegraphics[scale=0.7]{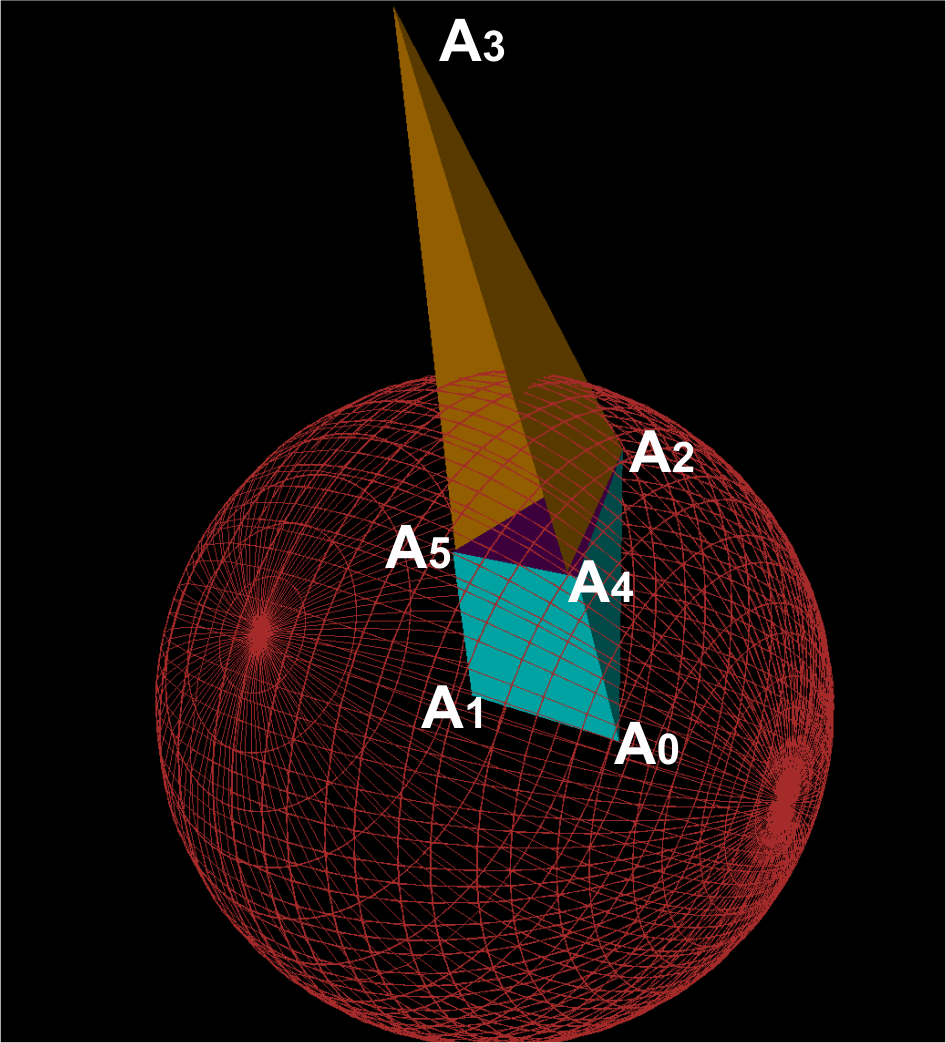}
            \caption{Truncated orthoscheme, where the truncating face $A_2 A_4 A_5$}
            \label{Truncated_Orthoscheme}
        \end{minipage}~~~~
        \begin{minipage}[b]{0.4\textwidth}
            \centering
            \includegraphics[scale=0.7]{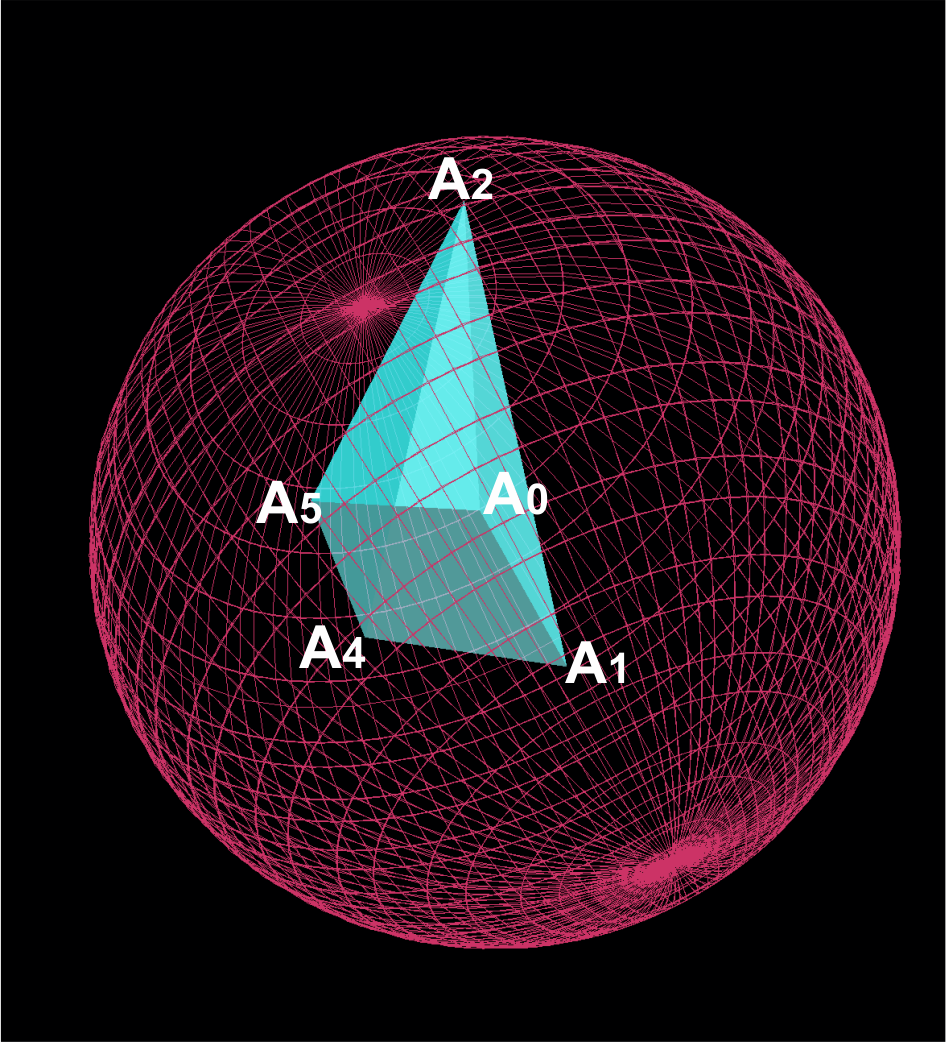}
            \caption{Truncated orthoscheme with the truncating face $A_2 A_4 A_5$}
            \label{Truncated_Orthoscheme_2}
        \end{minipage}
    \end{center}
\end{figure}

Our result is described in the Beltrami-Cayley-Klein-Sphere model of $3$ dimensional 
hyperbolic space $\HYP$, see Fig.\ref{Structure_Orthoscheme}, \ref{Orthoscheme},  \ref{Truncated_Orthoscheme}, and \ref{Truncated_Orthoscheme_2}. 
In Fig.\ref{Structure_Orthoscheme}, we provide the basic structure of a truncated orthoscheme with two pairs of parallel faces that 
intersect each other at the infinity. These intersection is described by two tangent lines of the model sphere, $k$ and $l$.
The corresponding Coxeter group is generated by reflections with respect
to its faces $h_0(\Bu_0), ~ h_1(\Bu_1),~ h_2(\Bu_2),~ h_3(\Bu_3),~ h_4(\Bu_4)$ (admitting the given Gram matrix $G$).
This group can be generated by reflections $\{T_0, T_1, T_2, T_3, T_4\}$, where $T_i$ $(i=0\dots4)$ denote the reflection with respect to the hyperplane $h_i$.
The reflections can be described in this model by the following formula:
\begin{equation*}
    T_i (\bx)=\bx-2 \langle \bx, \bu_i \rangle  \bu_i, ~~~~\text{where}~ i=0,1,\dots,4.
\end{equation*}
The computer visualization of the truncated orthoschemes are given in Fig.\ref{Orthoscheme}, \ref{Truncated_Orthoscheme}, and \ref{Truncated_Orthoscheme_2}.
\section{On sphere packings}
\begin{figure}[h!]
    \begin{center}
        \begin{minipage}[b]{0.4\textwidth}
            \centering
            \includegraphics[scale=0.27]{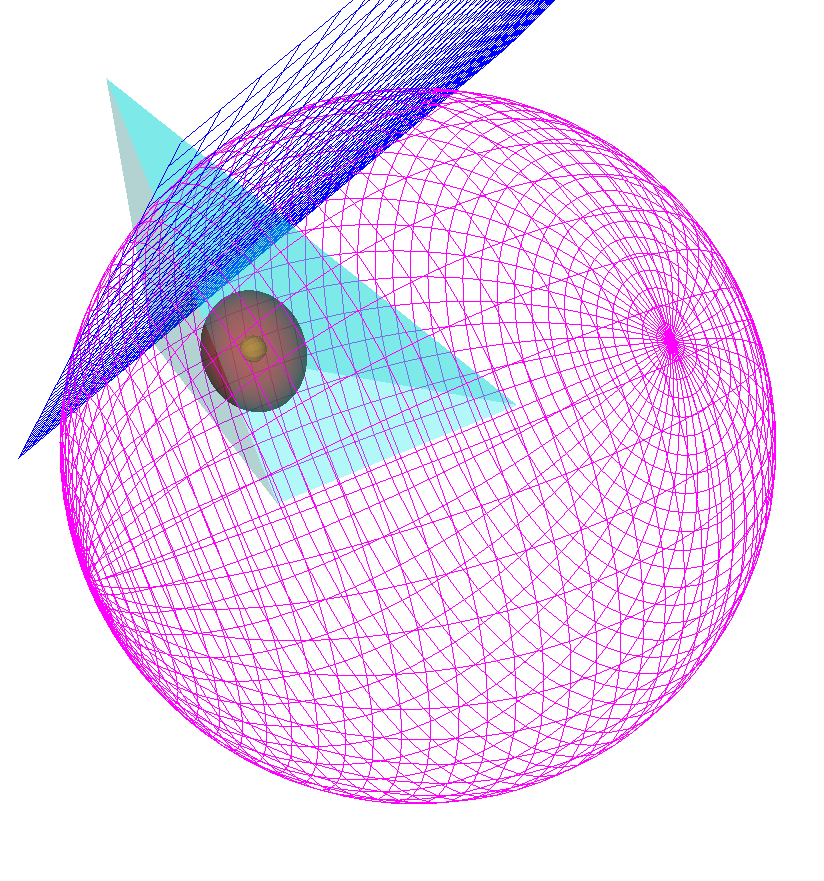}
            \caption{Optimum inball in the truncated orthoscheme $\{ \infty, 3, 3, \infty \}$}
            \label{inball}
        \end{minipage}~~~~~~
        \begin{minipage}[b]{0.4\textwidth}
            \centering
            \includegraphics[scale=0.27]{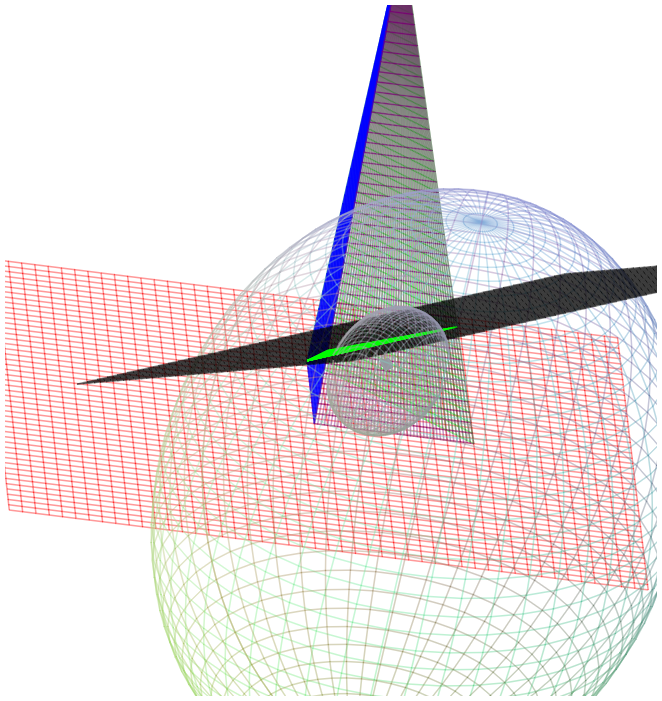}
            \caption{The inball intersects the truncating hyperplane $h_4$}
            \label{truncated_inball}
        \end{minipage}
    \end{center}
\end{figure}

In constructing the inball, the largest inscribed classical ball, in the truncated orthoscheme, we follow the procedure of \cite{J14}, 
that constructed the inscribed sphere into a polyhedra in hyperbolic spaces $\mathbb{H}^n$. 
For $i,j \in \{ 0, 1, 2, 3 \},~i\neq j$, let $S_{ij}(\bs_{ij})$ be hyperbolic hyperplane given by the following formula where $\bu_i$ are the above unit pols (normal vectors) of 
forms $\Bu_i$: 
\begin{equation}
    \Bs_{ij}:=(\bu_i-\bu_j)^{\perp}.
\end{equation}
Basically, the hyperplane $S_{ij}(\bs_{ij})$ is a bisector hyperplane that divide the dihedral angle between $h_i$ and $h_j$ into two equal parts.
Based on this notion, we construct a set of vectors as follows
\begin{equation*}
    \bs_{i}=\bu_{i}-\bu_{i+1},
\end{equation*}
for $i=0,1,2$. Note that the forms $\bs_{0}, \bs_{1}, \bs_{2}$ are linearly independent in $\mathbb{R}^{4}$. Now, we have three hyperbolic bisector hyperplanes $S_i(\bs_i)$,
where $i=0,1,2$.
Let $\hat{X}(\hat{\bx})$ be the common point of the bisector hyperplanes:
\begin{equation}
    \hat{\bx}:=\bigcap_{i=0,1,2}^{}S_i{(\bs_i)}.
\end{equation}
That means, we consider the equations system 
\begin{equation} \label{system to center}
    \langle \hat{\bx}, \bs_i \rangle =0, ~~~i=0,1,2.
\end{equation}
We consider that solution where $\hat{\bx}$ is a proper point in the model and lies in the given orthoscheme. 
The visualization of optimum inball in truncated orthoscheme $\{ \infty, 3, 3, \infty \}$ is given on Fig.\ref{inball}.
The problem might be occurred when the inball intersect the truncating hyperplane $\Bs_4$ (see Fig \ref{truncated_inball}). 
In that situation, the inradius is greater than the distance between the center of the inball and the truncating hyperplane $h_4$, i.e 
\begin{align*}
    r&=d_{\mathcal{H}}(X,h_0)=d_{\mathcal{H}}(X,h_1)=d_{\mathcal{H}}(X,h_2)=d_{\mathcal{H}}(X,h_3)\\
    &\ge d_{\mathcal{H}}(X,h_4).
\end{align*}
The three bisectors can be re-chosen in such that the inball does not intersect the omitting hyperplane.  There are $5$ complete combinations of bisectors for constructing the candidate of inball center. 
The complete packings denisties of inball packing (and their optimum density) can be found in \cite{YSz}, 
that gave the optimum packing density $\approx 0.2623649$, attained by sphere packing in $(\infty,3,3,\infty)$.
\section{On horosphere packings}
A horoshere in the hyperbolic geometry is the surface orthogonal to
the set of parallel lines, passing through the same point on the absolute quadratic surface $\partial \mathbb{H}^n$ (at present $n=3$)
(simply absolute) of the hyperbolic space $\HYP$.

We introduce an usual projective coordinate system using vector
basis $\boldsymbol{b}_i \ (i=0,1,2,3)$ for $\mathcal{P}^3$ where the
coordinates of center of the model is $C=(1,0,0,0)$. We pick an
arbitrary point at infinity to be $A_3=(1,0,0,1)$.

As it is known, the equation of a horosphere with center
$A_3=(1,0,0,1)$ through point $S=(1,0,0,s)$ $(s \in (-1,1))$ is 
\begin{gather}
\frac{{(s-1)}^2}{1-s^2}(-x^0 x^0 +x^1 x^1+x^2 x^2+x^3 x^3)+{(x^0-x^3)}^2=0 \notag
\end{gather}
This surface can be described in given usual Cartesian coordinate system by the following formula  
\begin{equation}\label{Key}
\displaystyle \dfrac{2(x^2+y^2)}{1-s} + \frac{4(z-(\frac{s+1}{2}))^2}{(1-s)^2}=1,
\end{equation}
where $x=\frac{x^1}{x^0},~y=\frac{x^2}{x^0},~z=\frac{x^3}{x^0}$. 

In computer visualization, it is very powerful to convert the horosphere equation into a polar coordinate system. We use the following conversion
\begin{align}\label{horopolareq1}
     x=\sqrt{\frac{1-s}{2}} \cos \theta \sin \phi, ~ &y=\sqrt{\frac{1-s}{2}} \sin \theta \sin \phi, ~z=\frac{1+s}{2}+&\frac{1-s}{2} \cos \phi,
\end{align}
where parameters $\theta \in [0,2 \pi)$, $\phi \in [0,\pi]$.

We will frequently be dealing with the arbitrary coordinate of orthoscheme (based on the set of unit normal Napier cycles). 
It would be useful if we place our arbitrary ideal point $A_2$ to $(1,0,0,1)$, by the following ideas.
We find the transformation $T$ (depend on two parameters) such that $T(A_2)=(1,0,0,1)$. Basically, this transformation 
$T$ is represented by a composition two rotations, i.e its representation matrix is described in the form
\begin{align*}
    \begin{bmatrix}1 & 0 & 0 & 0\\0 & \cos{\left(\phi \right)} & - \sin{\left(\phi \right)} & 0\\0 & \sin{\left(\phi \right)} & \cos{\left(\phi \right)} 
    & 0\\0 & 0 & 0 & 1\end{bmatrix} \left[\begin{matrix}1 & 0 & 0 & 0\\0 & 1 & 0 & 0\\0 & 0 & \cos{\left(\theta \right)} &  -\sin{\left(\theta \right)}\\
    0 & 0 & \sin{\left(\theta \right)} & \cos{\left(\theta \right)}\end{matrix}\right].
\end{align*}
Having the coordinate of ideal points $A_2$, one could determine the value of parameters $\phi$, and $\theta$.
Finally, by using transformation $T$, we move our truncated orthoscheme such that the ideal vertex $A_2$ has the coordinates $(1,0,0,1)$. 
Then one could apply the formulas (\ref{Key}), and (\ref{horopolareq1}) to construct the horosphere equations. \\

We separate our discussion into two cases depending on the number of vertices lying at infinity. 
First, we consider that case if the truncated orthoschemes have only one ideal vertex i.e 
$A_2$. In this case, we attach a horosphere centred at ideal vertex, $A_2$.\\~\\
While on the second case we shall consider the situation where two vertices of truncated orthoschemes lie on the absolute 
they are happened in cases $\{\infty,3,6, \infty\}$, $\{\infty,4,4, \infty\}$, and $\{\infty,6,3, \infty\}$. 
In each case, we could put a horosphere centered at either ideal vertex. We can also attach two horospheres altogether, where they are touching each other at edge $A_0 A_2$.\\
\subsubsection{Packings with one horosphere}
We have some truncated orthoschemes given with Schl\"afli symbols such that it has only one point at the infinity, 
e.g : $\{\infty,3,3,\infty\}$, $\{\infty,3,4, \infty\}$, $\{\infty,3,5, \infty\}$, $\{\infty,4,3, \infty\}$, $\{\infty,5,3, \infty\}$. 
However, if the truncated orthoscheme has two ideal vertices of truncated orthoscheme we can also study the corresponding horosphere packing centered at one either of these vertices. 

It is clear that the densest horoshpere packing configuration would be reached whenever this horosphere (horoball) with center $A_2$ touch the opposite face
(represented by hyperplane $h_2$). One could simply take the projection of $A_2$ onto its opposite face by the projection formula
$P_{h_2}(\ba_2)=\ba_2-\langle \ba_2, \bu_2 \rangle \bu_2 $.\\
The optimal horosphere should contains the point $P_{h_2}(A_2)$ therefore we can determine the parameter $s$ and so the equation (\ref{Key}) of the horosphere.

\begin{figure}[h!]
    \begin{center}
        \begin{minipage}[b]{0.4\textwidth}
            \includegraphics[scale=0.36]{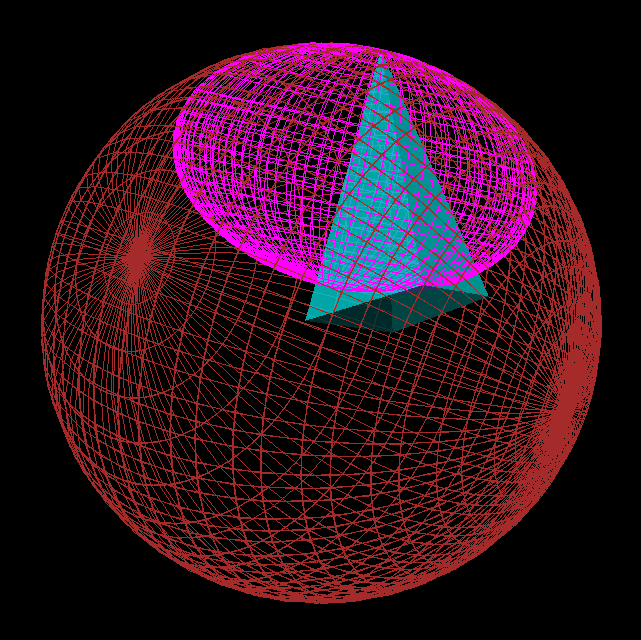}
            \vspace{0.5 cm}
            \caption{The largest horoball related to truncated orthoscheme of tiling $\{\infty,3,3,\infty\}$}
            \label{horo33_1}
        \end{minipage}~~~~~~~~~
        \begin{minipage}[b]{0.4\textwidth}
            \includegraphics[scale=0.37]{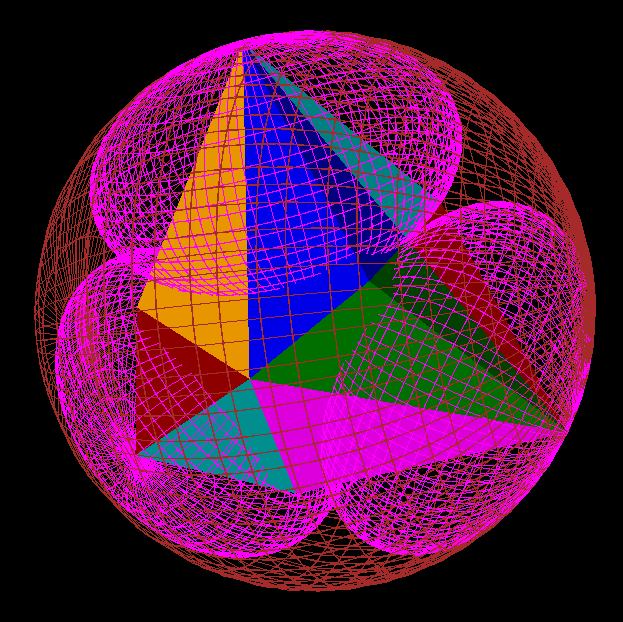}
            \caption{The neighbouring orthoschemes and horospheres configuration (first crown) related to truncated orthoscheme tiling $\{\infty,3,3,\infty\}$ }
            \label{horo33_2}
        \end{minipage}
    \end{center}
\end{figure}

\begin{figure}[h!]
    \begin{center}
        \begin{minipage}[b]{0.5\textwidth}
            \centering
            \includegraphics[scale=0.4]{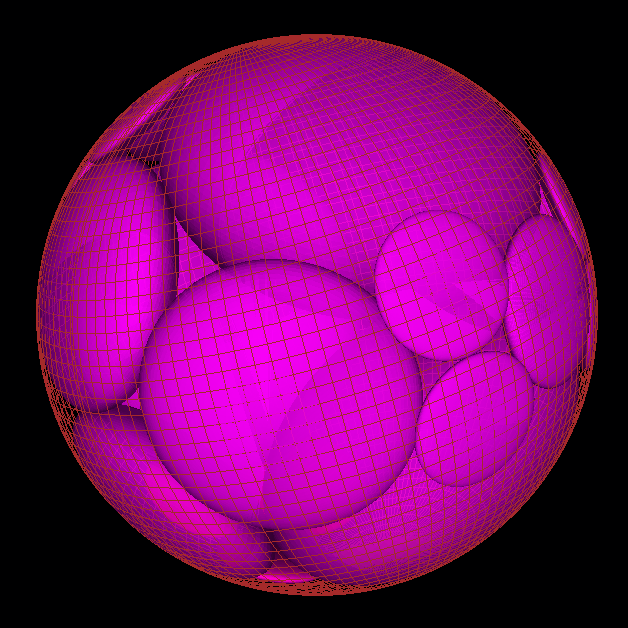}
            \caption{The first-third crowns of neighbouring horosphere configurations related to truncated orthoscheme tiling $\{\infty,3,3,\infty\}$}
            \label{horo33_3}
        \end{minipage}
    \end{center}
\end{figure}

We provide the computer visualization of optimum horospheres packing, attained by truncated orthoscheme tilings 
$\{\infty, 3,3, \infty\}$, in Fig.~\ref{horo33_1}, Fig.~\ref{horo33_2}, and Fig.~\ref{horo33_3}. The optimum packing density is $ \approx 0.8188080$, see \cite{YSz}.
\subsubsection{Packing with two horospheres}
Now, we focus on the orthoscheme tiling with the Schl\"{a}fli symbol e.g $\{\infty, 3, 6, \infty\}$, $\{\infty,4,4, \infty\}$, and $\{\infty, 6, 3, \infty\}$.
\begin{figure}[h!]
    \begin{center}
        \begin{minipage}[b]{0.5\textwidth}
            \centering
            \includegraphics[scale=0.7]{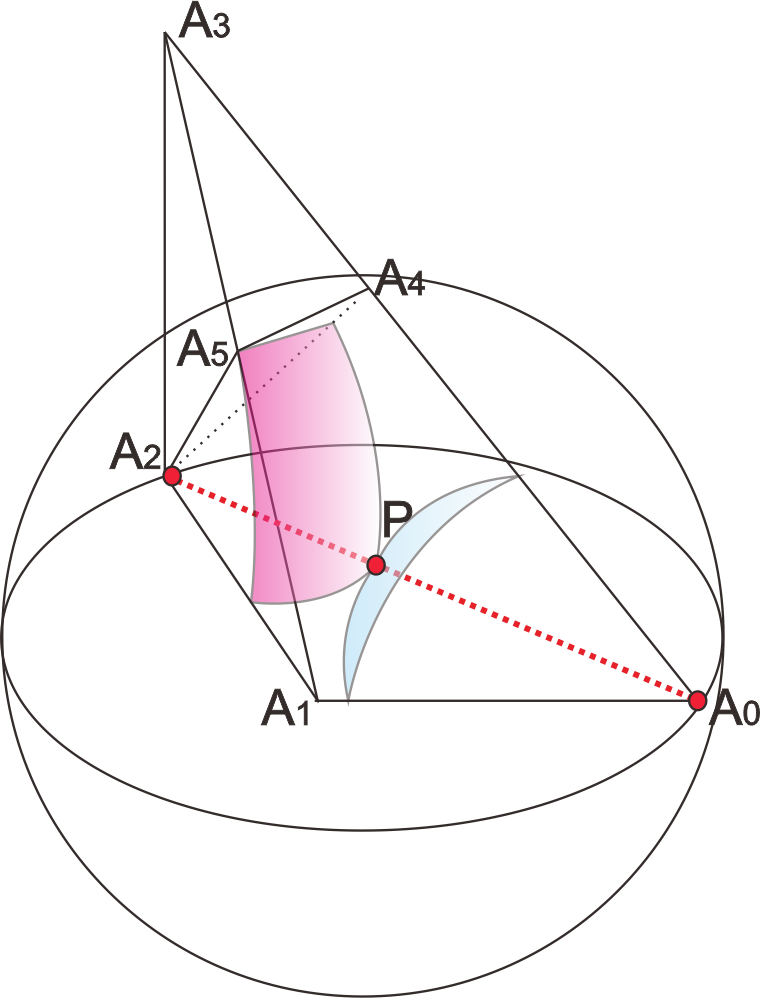}
            \caption{Two horospheres, $\mathcal{B}_0$ and  $\mathcal{B}_2$, that touch each other at a point lying on the edge $A_0A_2$}
            \label{double}
        \end{minipage}
   \end{center}
\end{figure}
In these cases given by the above Schl\"{a}fli symbols, we consider the generating truncated orthoscheme with $2$ ideal vertices 
i.e $A_2$ and $A_0$ lie at the infinity. We construct two horospheres centered at $A_2$ and $A_0$. We denote the horospheres centered at $A_2$ and $A_0$ 
with $\mathcal{B}_2$ and $\mathcal{B}_0$, respectively. It is clear, that in the optimal case these horospheres would touch each other at the edge joining 
$A_2$ and $A_0$, see Fig.(\ref{double}). 
The question is at which point of the edge the optimal (densest packing) configuration occurs if we move the touching point along edge joining $A_2$ and $A_0$. \\
\begin{rem}: 
These two horospheres could not intersect over their opposite faces, therefore there will be a restriction for the movement of the touching point.
\end{rem}
We can parameterize the possible movement of the touching point $P$, see Fig.~\ref{double}, e.g. $P(\b{r}(t))=(1-t)\ba_2+t\cdot \ba_0$. 
Then, for every possible $t$, we have parameters $s_i$ ($i=1,2)$ related to the both horospheres.\\

\vspace{3mm}

\textbf{Optimal horoball packing of tiling $\{\infty,4,4,\infty\}$} \\

\vspace{3mm}

In this situation, we have quite interesting structure,
we obtain that the possible values of $t$ are $t\in [\approx 0.2150 < t < \approx 0.3497]$.
We can compute the volumes of horoball sectors as the functions of $t$. It is just similar to the previous case, 
the volume function of horoball sectors centered at $A_2$ is a monotonic increasing function 
of $t$ if the touching point moving with direction $A_2$ to $A_0$ while the volume function of horoball sectors centered at $A_0$ is decreasing in this situation.
 
\begin{figure}[h!]
    \begin{center}
        \begin{minipage}[b]{0.5\textwidth}
            \centering
            \includegraphics[scale=0.6]{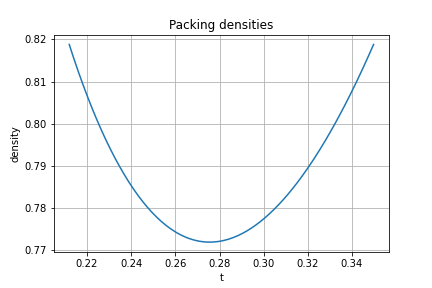}
            \caption{The plot of density function for all possible $t$ in case $\{\infty,4,4,\infty\}$}
            \label{dens44}
        \end{minipage}
    \end{center}
\end{figure}

In this case, we proved (see \cite{YSz}) that the density is increasing as a function of $t$, see Fig.~\ref{dens44}. 
Furthermore, the maximum density is attained when $t$ is largest, i.e when the horosphere centered at $A_2$ touches the opposite face.
\begin{equation*}
\delta_{opt}(\mathcal{B}({4,4,t})=\sup \delta(\mathcal{B}({4,4,t})) \approx 0.8188081.
\end{equation*}

\vspace{3mm}

\textbf{Optimal horoball packing of tilings $\{\infty,3,6,\infty\}$ and $\{\infty,6,3,\infty\}$}\\

\vspace{3mm}

We visualize similarly to the above case the densest horosphere (horoball) packings related to the truncated orthoscheme 
tilings with Schl\"afli symbol $\{\infty,3,6,\infty\}$ and $\{\infty,6,3,\infty\}$.
\begin{figure}[h!]
    \begin{center}
        \begin{minipage}[b]{0.4\textwidth}
            \centering
            \includegraphics[scale=0.7]{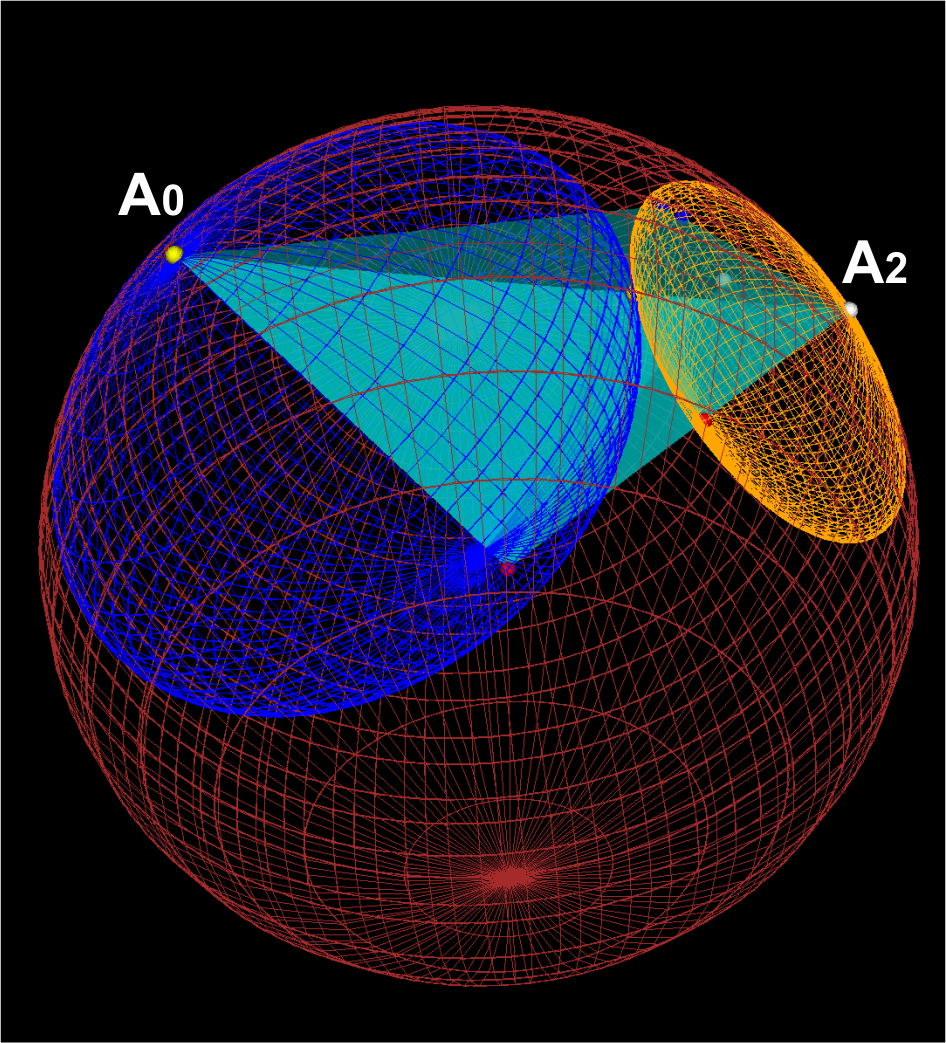}
            \vspace{1 pt}
            \caption{Two horospheres, $\mathcal{B}_0$ and  $\mathcal{B}_2$, that touch each other at the point lying on $A_0A_2$ related to tiling $\{\infty,3,6,\infty\}$.}
            \label{h36_0}
        \end{minipage}~~~~~~~
        \begin{minipage}[b]{0.4\textwidth}
            \centering
            \includegraphics[scale=0.35]{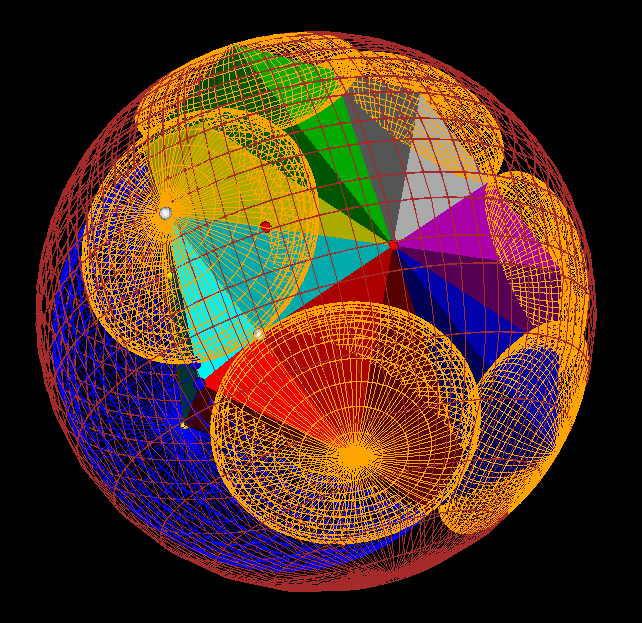}
            \caption{Adjacent orthoschemes and the corresponding horosphere configuration (first crown) related to truncated orthoscheme tiling$\{\infty,3,6,\infty\}$}
            \label{h36_1}
        \end{minipage}
    \end{center}
\end{figure}
\begin{figure}[h!]
    \begin{center}
        \begin{minipage}[b]{0.4\textwidth}
            \centering
            \includegraphics[scale=0.355]{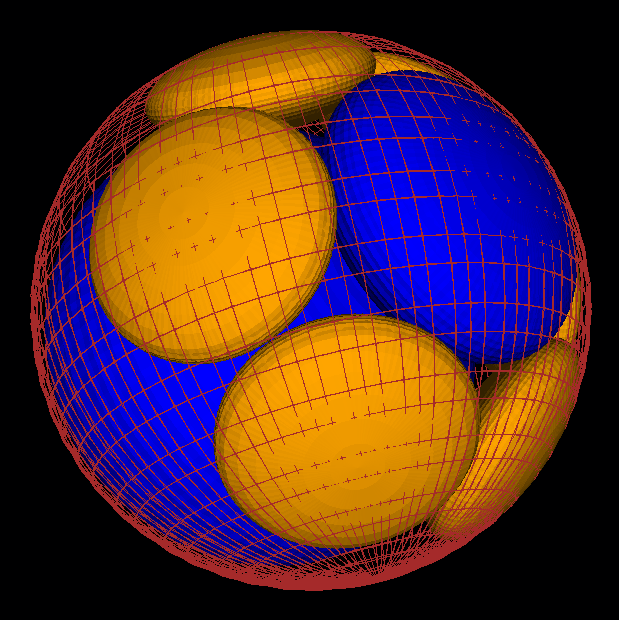}
            \vspace{1 pt}
            \caption{The horosphere configuration (first crown) related to truncated orthoscheme tiling $\{\infty,3,6,\infty\}$}
            \label{h36_2}
        \end{minipage}~~~~~
        \begin{minipage}[b]{0.4\textwidth}
            \centering
            \includegraphics[scale=0.30]{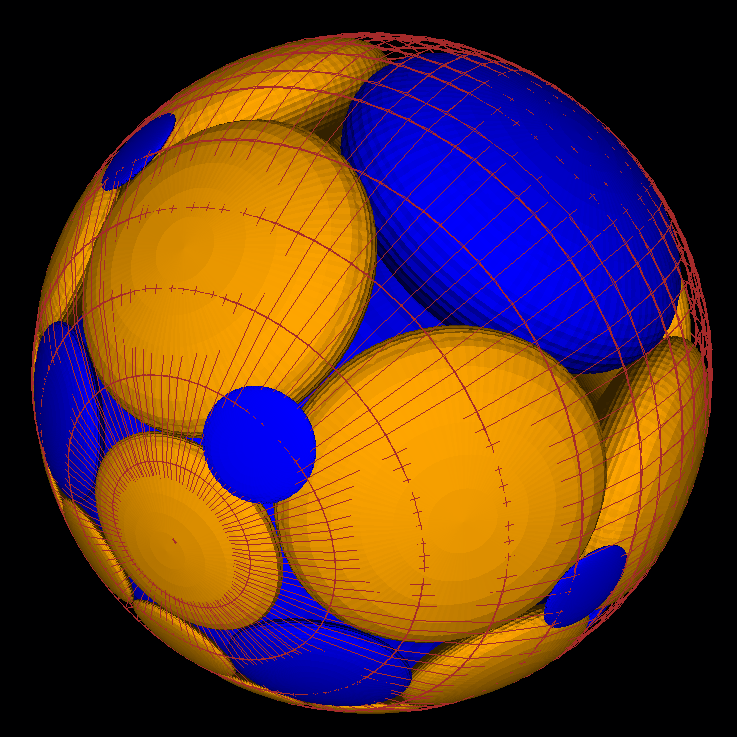}
            \caption{The optimum packing density horospheres configuration (first-third crown) related to the orthoscheme tiling $\{\infty,3,6,\infty\}$.}
            \label{h36_3}
        \end{minipage}
    \end{center}
\end{figure}
There are some basic facts occurred in the orthoschemes with Schl\"afli symbol $\{\infty,3,6,\infty\}$, $\{\infty,6,3,\infty\}$.
\begin{enumerate}
    \item In this situation, there is in each case only one possible value of parameters $t$, $t_{(3,6)}\approx 0.2119416,~ t_{(6,3)}$, $\approx 0.5745582.$
    \item If $(q,r)=(3,6)$ then the optimal horosphere $\mathcal{B}_2$ touches the plane $h_2$ and $\mathcal{B}_0$ touches the face $h_0$ and if 
    $(q,r)=(6,3)$ $\mathcal{B}_2$ touches the plane $h_2$ and $\mathcal{B}_0$ touches the face $h_4$.
    \item The packing density of these two configurations are the same, $\approx 0.8413392$, see \cite{YSz}.
\end{enumerate}
Finally, we give the computer visualization in Fig.~\ref{h36_0}-\ref{h36_2} related to Coxeter tiling $\{\infty,3,6,\infty\}$ and in Fig.~\ref{h63_0}-\ref{h63_3} for Coxeter 
tiling $\{\infty,6,3,\infty\}$.

In our opinion, non-Euclidean tilings and packings and their investigations will play an important role in the research of
material structure in the near future. Therefore, we also consider it important to visualize them in order to get to know them. 
\begin{figure}[h!]
    \begin{center}
        \begin{minipage}[b]{0.4\textwidth}
            \centering
            \includegraphics[scale=0.287]{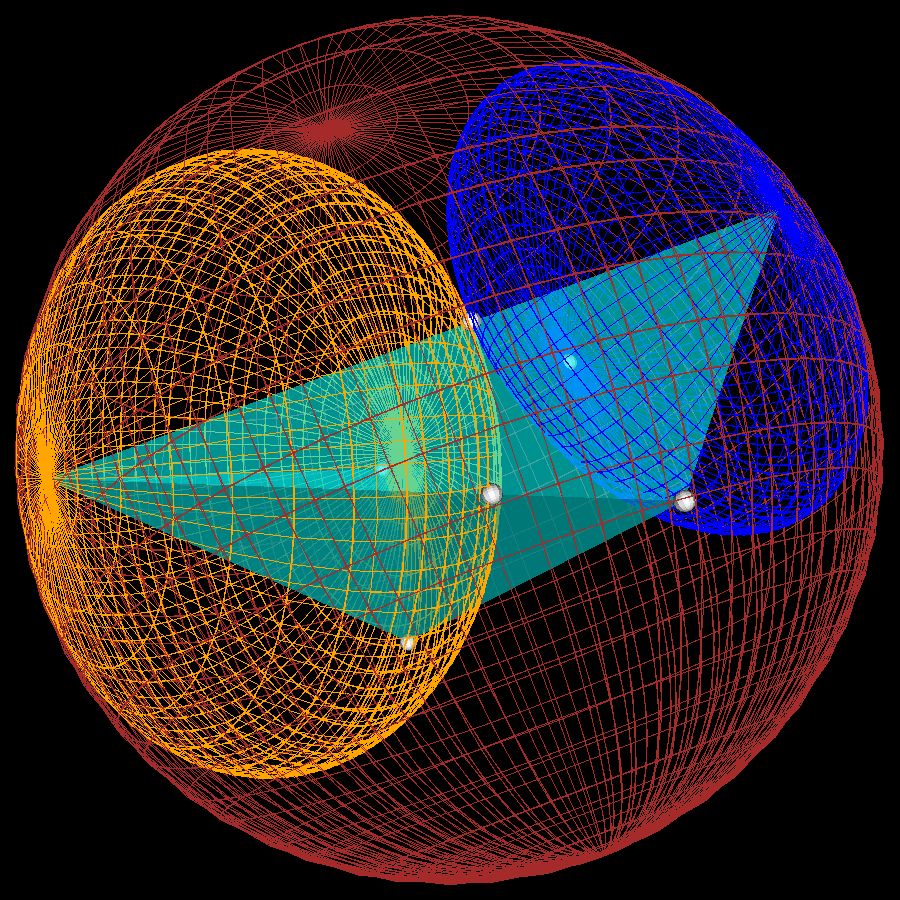}
            \caption{Two horospheres, $\mathcal{B}_0$ and  $\mathcal{B}_2$, that touch each other at a point lies on edge $A_0A_2$ related to tiling $\{\infty,6,3,\infty\}$.}
            \label{h63_0}
        \end{minipage}~~~~~~~~~~~~~~
        \begin{minipage}[b]{0.4\textwidth}
            \centering
            \includegraphics[scale=0.4]{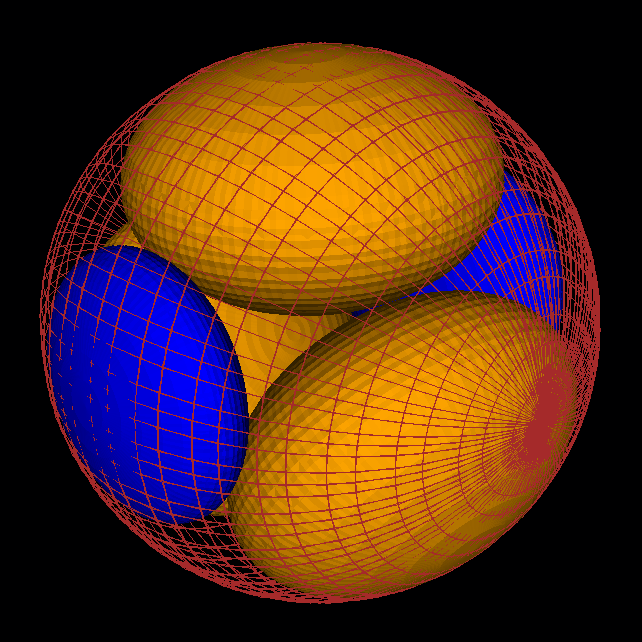}
            \caption{The horosphere configuration (first crown) related to truncated orthoscheme tiling $\{\infty,6,3,\infty\}$}
            \label{h63_1}
        \end{minipage}
    \end{center}
\end{figure}
\begin{figure}[h!]
    \begin{center}
        \begin{minipage}[b]{0.4\textwidth}
            \centering
            \includegraphics[scale=0.4]{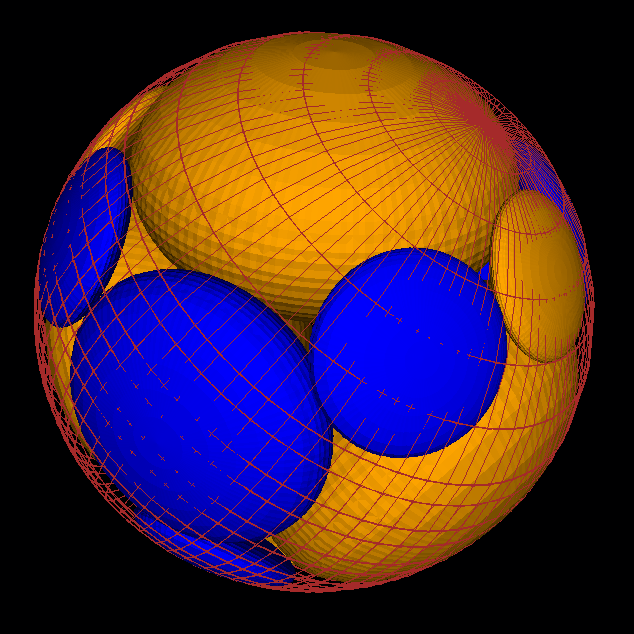}
            \caption{The optimum packing density horosphere configuration (first-second crown) related to the orthoscheme tiling $\{\infty,6,3,\infty\}$.}
            \label{h63_2}
        \end{minipage}~~~~~~~~~~~~~~
        \begin{minipage}[b]{0.4\textwidth}
            \centering
            \includegraphics[scale=0.4]{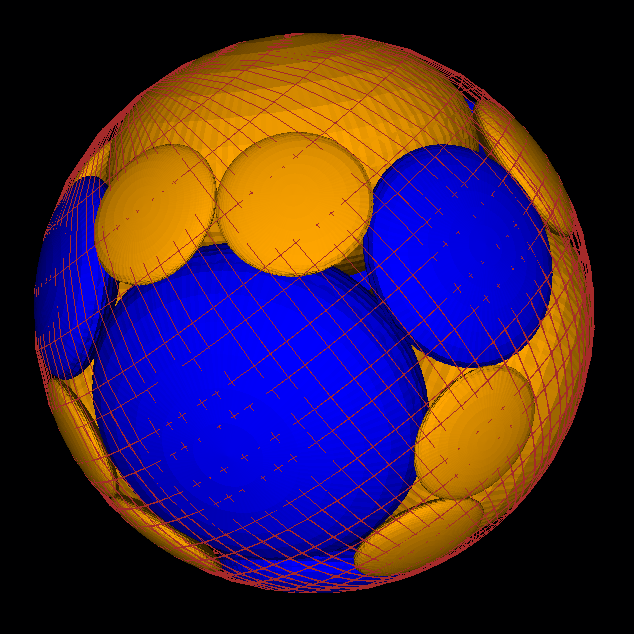}
            \caption{The optimum packing density horosphere configuration (first-third crown) related to the orthoscheme tiling $\{\infty,6,3,\infty\}$.}
            \label{h63_3}
        \end{minipage}
    \end{center}
\end{figure}
%

\end{document}